\documentclass{amsart}

\newtheorem{theorem}{Theorem}[section]
\newtheorem{lemma}[theorem]{Lemma}
\newtheorem{corollary}[theorem]{Corollary}

\theoremstyle{definition}

\newtheorem{example}[theorem]{Example}

\theoremstyle{remark}
\newtheorem{remark}[theorem]{Remark}

\numberwithin{equation}{section}

\begin{document}

\title[Nevanlinna characteristic of  ${f(z+\eta)}$ and
difference equations] {On the Nevanlinna characteristic of
$f(z+\eta)$ and difference equations in the complex plane}

\author{Yik-Man Chiang}
\address{Department of Mathematics, Hong Kong University of Science and
Technology, Clear Water Bay, Kowloon, Hong Kong,\ P. R. China.}
\email{machiang@ust.hk}
\thanks{This research was supported in part by the Research Grants Council
of the Hong Kong Special Administrative Region, China
(HKUST6135/01P)}

\author{Shao-Ji Feng}
\address{Academy of Mathematics and Systems Science, Chinese Academy
of Sciences, Beijing, 100080, P. R. China.}
%\curraddr{Department of Mathematics, Hong Kong University of
%Science and Technology, Clear Water Bay, Kowloon, Hong Kong,\ P. R. China.}
\email{fsj@amss.ac.cn}
\thanks{The second author was also partially supported by National Natural Science
Foundation of China (Grant No. 10501044) and by the HKUST PDF
Matching Fund.}

%    General info

\date{To appear in The Ramanujan J., $14^{\textrm{th}}$ Nov. 2005.}

%\dedicatory{This paper is dedicated to our authors.}


\begin{abstract}
We investigate the growth of the Nevanlinna Characteristic of
$f(z+\eta)$ for a fixed $\eta\in\textbf{C}$ in this paper. In
particular, we obtain a precise asymptotic relation between
$T\big(r,\, f(z+\eta)\big)$ and $T(r,\,f)$, which is only true for
finite order meromorphic functions. We have also obtained the
proximity function and pointwise estimates of $f(z+\eta)/f(z)$
which is a discrete version of the classical logarithmic
derivative estimates of $f(z)$. We apply these results to give new
growth estimates of meromorphic solutions to higher order linear
difference equations. This also allows us to solve an old problem
of Whittaker \cite{Whit} concerning a first order difference
equation. We show by giving a number of examples that all of our
results are best possible in certain senses. Finally, we give a
direct proof of a result in Ablowitz, Halburd and Herbst
\cite{AHH} concerning integrable difference equations.
\end{abstract}

\maketitle

%\section*{This is an unnumbered first-level section head}
%This is an example of an unnumbered first-level heading.

%% The correct journal style for \specialsection is all uppercase; a known bug
%% in amsart.cls prevents this, so input must be uppercase until it is fixed.
%\specialsection*{This is a Special Section Head}
%\specialsection*{THIS IS A SPECIAL SECTION HEAD}
%This is an example of a special section head%
%%%%%%%%%%%%%%%%%%%%%%%%%%%%%%%%%%%%%%%%%%%%%%%%%%%%%%%%%%%%%%%%%%%%%%%%
%\footnote{Here is an example of a footnote. Notice that this footnote
%text is running on so that it can stand as an example of how a footnote
%with separate paragraphs should be written.
%\par
%And here is the beginning of the second paragraph.}%
%%%%%%%%%%%%%%%%%%%%%%%%%%%%%%%%%%%%%%%%%%%%%%%%%%%%%%%%%%%%%%%%%%%%%%%%

\section{Introduction}
A function $f(z)$ is called meromorphic if it is analytic in the
complex plane \textbf{C} except at isolated poles. In what
follows, we assume the reader is familiar with the basic notion of
Nevanlinna's value distribution theory (see e.g. \cite{Nev},
\cite{Hay64}).

 Recently,
there has been renewed interests in difference (discrete)
equations in the complex plane $\textbf{C}$ (\cite{AHH},
\cite{CR}, \cite{Hay90}, \cite{HKLRT}, \cite{Is-Ya}, \cite{Ruij};
see also \cite{BL}). In particular, and most noticeably, is the
proposal by Ablowitz, Halburd and Herbst \cite{AHH} to use the
notion of order of growth of meromorphic functions in the sense of
classical Nevanlinna theory \cite{Hay64} as a detector of
\textit{integrability} (i.e., solvability) of second order
non-linear difference equations in \textbf{C}. In particular, they
showed in \cite{AHH} that if the difference equation
\begin{equation}
\label{E:Malmquist}
    f(z+1)+f(z-1)=R(z,\,f(z))=\frac{a_0(z)+a_1(z)f(z)+\cdots+a_p(z)f(z)^p}
    {b_0(z)+b_1(z)f(z)+\cdots+b_q(z)f(z)^q}
\end{equation}
admits a finite order meromorphic solution, then $\max(p,\, q)\le
2$. It is proposed in \cite{AHH} that a difference equation admits
a \textit{finite order} meromorphic solution is a strong
indication of \textit{integrability} of the equation. It is known
that when $\max(p,\, q)\le 2$, the equation (\ref{E:Malmquist})
includes the well-known discrete Painlev\'e equations which are
prime examples of integrable second order difference equations
(\cite{AHH}, \cite{GLS}). More discussion about the integrability
of difference equations will be relegated to \S10.

In contrast to differential equations, non-linear difference
equations often admit global meromorphic solutions (\cite{Shim},
\cite{Yana1980}) and hence Nevanlinna's value distribution theory
is applicable. The classical growth comparisons between $f(z)$ and
$f^\prime(z)$ have important applications to differential
equations (see e.g.  \cite{Gund-2}, \cite{Hay96}). In the case of
applying Nevanlinna theory to difference equations, one of the
most basic questions is the growth comparison between
$T(r,\,f(z+1))$ and $T(r,\,f(z))$. It is shown in \cite[p.
66]{GO}, that for an arbitrary $b\not=0$, the following
inequalities\footnote{The inequalities (\ref{E:Goldberg}) can be
easily derived from corresponding inequalities for the
Ahlfors-Shimizu characteristic function \cite{GO}.}
\begin{equation}
\label{E:Goldberg}
    \big(1+o(1)\big)T\big(r-|b|,\, f(z)\big)\le T\big(r,\,
    f(z+b)\big)\le  \big(1+o(1)\big)T\big(r+|b|,\, f(z)\big)
\end{equation}
hold as $r\to \infty$ for a general meromorphic function. Let
$\eta$ be a non-zero complex number, we shall prove the following
precise asymptotic relation
\begin{equation}
\label{E:char-asy}
    T(r,\,f(z))\sim T\big(r,\, f(z+\eta)\big),
\footnote{If $\phi(r)$ and $\psi(r)$ are two real functions of a
real variable $r$, then $\phi\sim\psi$ means $\phi(r)/\psi(r)$
tends to unity as $r\to+\infty$, see \cite[p. 4]{Olv}}
\end{equation}
 holds for finite
order meromorphic functions.

 In this paper, instead of
 considering non-linear difference equations, we shall, however, concentrate
 ourselves on the value distribution
properties of $f(z+\eta)$ and related expressions and their
applications to \textit{linear difference equations}. It turns out
that the results we obtained not only allow us to give new results
on linear difference equations, but they also allow us to give a
direct proof of the Nevanlinna-type theorems in \cite{AHH},
including that for the equation (\ref{E:Malmquist}).

Although there has been considerable progress of the knowledge on
the growth of meromorphic solutions to $q$-difference equations
(e.g. \cite{BH}, \cite{BIY-1}, \cite{BIY-2}, \cite{HLRY},
\cite{Ram}), relatively little is known for the growth of
meromorphic solutions to even the first order difference equation
\begin{equation}
\label{E:1st-order}
    F(z+1)=\Psi(z) F(z),
\end{equation}
where $\Psi(z)$ is a meromorphic coefficient. If $\Psi(z)$ is a
rational function, then a solution is given by
\begin{equation}
\label{E:gamma-generalized}
    F(z)=e^{az}\frac{\prod_{j=1}^{m} \Gamma(z-b_j)}{\prod_{k=1}^{n}
    \Gamma(z-c_k)}
\end{equation}
for some suitable choices of constants $a,\, b_j,\, c_k\in
\textbf{C}$. Thus the solution has order 1.

Given a finite order meromorphic coefficient $\Psi(z)$, Whittaker
\cite[\S 6]{Whit} explicitly constructed a meromorphic solution
$F(z)$ of order $\le \sigma(\Psi)+1$ to (\ref{E:1st-order}).  On
the other hand, let $\Pi (z)$ be a periodic meromorphic function
of period 1, then the product $\Pi(z)F(z)$ again satisfies the
equation (\ref{E:1st-order}). Thus, we can get at most a lower
bound order estimate for a general meromorphic solution to
(\ref{E:1st-order}). In this paper, we shall settle the lower
bound order estimate for the Whittaker problem and, moreover, for
the higher order linear difference equations
\begin{equation}
\label{E:higher-order-linear}
    A_n(z)f(z+n)+\cdots+A_{1}(z)f(z+1) +A_0f(z)=0
\end{equation}
with certain entire coefficients $A_j(z)$. It will be shown that
the order of growth of a meromorphic solution is one larger than
the order of the  dominant coefficient amongst the $A_j(z)$.
Examples are given to demonstrate that the lower bound we obtain
is the best possible.

It turns out that the fundamental estimate (\ref{E:char-asy})
needs both the
\begin{equation}
\label{E:discrete-counting}
    N\big(r,\, f(z+\eta)\big)\sim N(r,\, f)
\end{equation}
for finite order meromorphic functions, as well as a version of
\textit{discrete} analogue of the classical logarithmic derivative
to be discussed below.

 It is
well-known that the following logarithmic derivative estimate
\begin{equation}
\label{E:log-derivative}
    m\left(r,\,\frac{f^\prime(z)}{f(z)}\right)=
    O\big(\log T(r,\,
    f)\big)=S(r,\, f),
\end{equation}
holds outside a possible set of finite linear measure, where the
notation  $S(r,\, f)$ means that the expression is of
$o\big(T(r,\, f))$. It shows that the proximity function of the
logarithmic derivative of $f(z)$ grows much slower than the
Nevanlinna characteristic function of $f(z)$. The \textit{
logarithm derivative lemma}, as it is often called, has numerous
applications in complex differential equations \cite{Laine} and it
also plays a crucial role in proving the celebrated Nevanlinna
Second Fundamental theorem \cite{Hay64}, \cite{Laine}. It is
generally recognized that the estimate (\ref{E:log-derivative}) is
amongst the deepest results in the value distribution theory. One
can also find other applications of it in \cite{Hinkk},
\cite{Lang}.

Let $\eta$ be a fixed complex number and $f(z)$ a meromorphic
function, we ask under what assumption on $f(z)$ do we have the
following difference analogue of (\ref{E:log-derivative})
\begin{equation}
\label{E:discrete-quotient}
    m\bigg(r,\,\displaystyle{\frac{f(z+\eta)}{f(z)}}\bigg)=S^*(r,\,f)
    \ ?
\end{equation}
Here the notation $S^*(r,\, f)$ means that the left hand side of
(\ref{E:discrete-quotient}) is of slower growth than $T(r,\, f)$
in some sense.

We shall give an answer to the question
(\ref{E:discrete-quotient}). More specifically, we show that if
$f(z)$ is a meromorphic function of finite order $\sigma$, then we
have
\begin{equation}
\label{E:discrete-quot-est}
    m\bigg(r,\,\frac{f(z)}{f(z+\eta)}\bigg)+
    m\bigg(r,\,\frac{f(z+\eta)}{f(z)}\bigg)=
    O\big(r^{\sigma-1+\varepsilon}\big)
\end{equation}
for an arbitrary $\varepsilon>0$. This estimate holds without any
exceptional set. Hence we obtain (\ref{E:discrete-quotient}) when
we choose $\varepsilon$ in (\ref{E:discrete-quot-est}) to be
sufficiently small. It is not difficult to see that it is
impossible for (\ref{E:discrete-quotient}) to hold for an
arbitrary meromorphic function. In fact,
(\ref{E:discrete-quotient}) fails to hold even for the simple
entire function $f(z)=e^{e^z}$ and $\eta\not= 2\pi ik,\, (k=1,\,
2,\,3,\cdots)$. After this paper is completed we
learnt\footnote{Mid--April 2005} that R.~ G.~Halburd and R. J.
Korhonen \cite{HK-1} have also obtained an essentially same
estimate (\ref{E:discrete-quot-est}), and its interesting
applications in \cite{HK-2} and \cite{HK-3}.

Although our problem regarding (\ref{E:discrete-quotient}) is
somewhat weaker than the (\ref{E:log-derivative}), we show that it
is already sufficient for our applications and, more importantly,
we shall show by examples that both the upper bounds and the
finite order assumption are the best possible.

The idea of the proof of (\ref{E:discrete-quot-est}) relies on an
application of the Poisson--Jensen formula \cite{Hay64}. The
formula also allows us to obtain pointwise estimates  for
$|{f(z+\eta)}/{f(z)}|$.

We recall that Gundersen \cite[Cor. 1]{Gund} has given a precise
pointwise estimate for the logarithmic derivative for a
meromorphic function $f(z)$ of order $\sigma$ to be
\begin{equation}
\label{E:pointwise-der}
    \left|\frac{f^\prime(z)}{f(z)}\right|\le
    |z|^{\sigma-1+\varepsilon}
\end{equation}
for all $|z|$ sufficiently large and outside some \textit{small
exceptional sets}. Our estimates allow us to show, amongst others,
the new upper bound
\begin{equation}
\label{E:pointwise-est}
    \left|\frac{f(z+\eta)}{f(z)}\right|\le \exp\{r^{\sigma-1+\varepsilon}\}
\end{equation}
where $|z|=r$ is sufficiently large and $|z|$ is also outside some
\textit{small exceptional sets}. We shall apply these estimates to
obtain new growth estimates of entire solutions to equation
(\ref{E:higher-order-linear}) with \textit{polynomial}
coefficients.

This paper is organized as follows. The main results concerning
the growth of the Nevanlinna characteristic of $f(z+\eta)$ will be
stated in \S2; some preliminary lemmas are stated and proved in
\S3. The proof of the main theorems are given in \S4 to \S7. We
will consider pointwise estimates such as  (\ref{E:pointwise-est})
in \S8. The applications of the main results to difference
equations are given in \S9, followed by a discussion of the
relation of our results to integrable difference equations in
\S10.

%\subsection{This is a numbered second-level section head}
%This is an example of a numbered second-level heading.

%\subsection*{This is an unnumbered second-level section head}
%This is an example of an unnumbered second-level heading.

%\subsubsection{This is a numbered third-level section head}
%This is an example of a numbered third-level heading.

%\subsubsection*{This is an unnumbered third-level section head}
%This is an example of an unnumbered third-level heading.

\section{Main Results on Nevanlinna Characteristics}
When $f(z)$ has a finite order of growth, we shall improve the
inequalities (\ref{E:Goldberg}) to the following theorem.

\begin{theorem}
\label{T:char} Let $f(z)$ be a meromorphic function with order
$\sigma=\sigma(f),\,\sigma<+\infty$, and let $\eta$ be a fixed non
zero complex number, then for each $\varepsilon>0$, we have
\begin{equation}
\label{E:char-est}
    T\big(r,f(z+\eta)\big)=T(r,\,f)+O(r^{\sigma-1+\varepsilon})+O(\log r).
\end{equation}
\end{theorem}
\medskip

It is interesting to compare the inequalities (\ref{E:Goldberg})
and (\ref{E:char-est}) with the following estimate given by
Ablowitz, Halburd and Herbst \cite[Lemma 1]{AHH} that given any
$\varepsilon>0$, then for all $r\ge 1/\varepsilon$, we have
\begin{equation}
\label{E:AHH-char-est}
    T\big(r, f(z\pm 1)\big)\le (1+\varepsilon) T\big(r+1,
    f(z)\big)+\kappa,
\end{equation}
where $\kappa$ is a constant. Thus our (\ref{E:char-est}) shows
that we have ``equality" in (\ref{E:AHH-char-est}) and that we can
choose $\varepsilon=0$ there, although we have a larger remainder
term in (\ref{E:char-est}). Although the technique used in
\cite{AHH} to obtain (\ref{E:AHH-char-est}) is different from that
of the (\ref{E:Goldberg}) in \cite{GO}, the latter has already
contained the (\ref{E:AHH-char-est}).

The above main theorem on Nevanlinna characteristic depends on the
following results.

\begin{theorem}
\label{T:counting} \ Let $f$ be a meromorphic function with
exponent of convergence of poles
$\lambda(\frac{1}{f})=\lambda<+\infty$, $\eta\not=0$ be fixed,
then for each $\varepsilon>0$,
\begin{equation}
\label{E:counting-est}
    N\big(r,f(z+\eta)\big)=N(r, f)+O(r^{\lambda-1+\varepsilon})+O(\log r).
\end{equation}
\end{theorem}
\medskip

The following example shows that the above Theorem is sharp in the
sense that (\ref{E:counting-est}) no longer hold for infinite
order meromorphic functions.

\begin{theorem}
\label{T:counting-example} There exists a meromorphic function
$f(z)$ of infinite order such that
\begin{equation}
\label{E:counting-example}
    \frac{N\big(r,\, f(z+1)\big)-N\big(r,\, f(z)\big)}{N\big(r,\,
    f(z)\big)}\ge 1
\end{equation}
as $r\to \infty$.
\end{theorem}
\medskip

\begin{theorem}
\label{T:discrete-quotient} Let $\alpha,\, R,\, R^\prime$ be real
numbers such that $0<\alpha< 1,\, 0<R$, and let $\eta$ be a
non-zero complex number. Then there is a positive constant
$C_\alpha$ depending only on $\alpha$ such that for a given
meromorphic function $f(z)$

we have, when $|z|=r,\, \max\{1,\, r+|\eta|\}<R<R^\prime$, the
estimate
\begin{equation}
\label{E:fund-est-I}
\begin{split}
    m\left(r,\,\frac{f(z+\eta)}{f(z)}\right)
    &+m\left(r,\,\frac{f(z)}{f(z+\eta)}\right)\\
    &\le\frac{2|\eta|R}{(R-r-|\eta|)^2}\left(m\big(R,\, f\big)+
    m\big(R,\, \frac{1}{f}\big)\right)+\\
    &\hskip5pt+\left(\frac{2R^\prime}{R^\prime-R}\right)\left(\frac{|\eta|}{R-r-|\eta|}
    +\frac{C_\alpha|\eta|^\alpha}{(1-\alpha)r^\alpha}\right)\left(N\big(R^\prime,\, f\big)+
    N\big(R^\prime,\, \frac{1}{f}\big)\right).
%+    &\hskip5pt+k\log(R/r),
\end{split}
\end{equation}
%where $k=|\kappa|+|\tau|$.
\end{theorem}
\medskip

We immediately deduce from (\ref{E:fund-est-I}) the following
corollary for finite order meromorphic functions.

\begin{corollary}
\label{C:discrete-quotient} Let $f(z)$ be a meromorphic function
of finite order $\sigma$ and let $\eta$ be a non-zero complex
number. Then for each $\varepsilon>0$, we have
\begin{equation}
\label{E:discrete-proximity}
    m\Big(r,\,\frac{f(z+\eta)}{f(z)}\Big)+m\Big(r,\,
    \frac{f(z)}{f(z+\eta)}\Big)=O(r^{\sigma-1+\varepsilon}).
\end{equation}
\end{corollary}
%\medskip

%\emph{Proof.}
\begin{proof}
Since $f(z)$ has finite order $\sigma(f)=\sigma<+\infty$, so
given $\varepsilon,\,0<\varepsilon<2$, we have
\[
T(r,\, f)=O(r^{\sigma+\frac{\varepsilon}{2}})
\]
for all $r$. We obtain (\ref{E:discrete-proximity}) by choosing
$\alpha=1-\frac{\varepsilon}{2}$, $R=2r,\, R^\prime=3r$ and
$r>\max\{|\eta|,\, 1/2\}$ in Theorem \ref{T:discrete-quotient}.
%We note that we cannot discard the $O(1)$ term in (\ref{E:discrete-proximity}) since the bound
%$O(r^{\sigma-1+\varepsilon})$ tends to zero when $0<\sigma<1$, and
%the last term in (\ref{E:fund-est-I}) becomes a constant after putting $R=2r$.
This completes the proof.
\end{proof}

We also deduce from (\ref{E:fund-est-I}) the following result.

\begin{corollary}
\label{T:discrete-quotient-2} Let $\eta_1,\, \eta_2$ be two
complex numbers such that $\eta_1\not=\eta_2$ and let $f(z)$ be a
finite order meromorphic function. Let $\sigma$ be the order of
$f(z)$, then for each $\varepsilon>0$, we have
\begin{equation}
\label{E:discrete-quotient-2}
    m\Big(r,\,
    \frac{f(z+\eta_1)}{f(z+\eta_2)}\Big)=O(r^{\sigma-1+\varepsilon}).
\end{equation}
\end{corollary}
\bigskip

We note that the estimate (\ref{E:discrete-proximity}) satisfies
the (\ref{E:discrete-quotient}) which is a discrete analogue of
(\ref{E:log-derivative}). This answers our question raised in
(\ref{E:discrete-quotient}) in the Introduction. We note the above
estimates do not hold when the order of $f(z)$ is infinite as
indicated in the following example. Hence they are the best
possible.

\begin{example}
\label{eg:one}
 Let $g(z)=e^{e^z},\, z=re^{i\theta}$.  We choose
$\eta$ to be real. It follows from \cite[p. 7]{Hay64} that
\begin{equation*}
%\begin{split}
    m\Big(r,\,\frac{g(z+\eta)}{g(z)}\Big)=(e^\eta-1)m(r,\, g)
    =(e^\eta-1)T(r,\, g)
    \sim (e^\eta-1)\frac{e^r}{\sqrt{2\pi^3r}}
%    \frac{1}{2\pi}\int_{-\pi}^{\pi}
%    e^{r\cos(\theta)}\phi(r\sin\theta)\,d\theta,
%\end{split}
\end{equation*}
as $r\to+\infty$.
\end{example}

The Example \ref{eg:one} demonstrates that $m(r,\,
f(z+\eta)/f(z))$ can grow as fast as the $m(r,\, f(z))$ itself for
an infinite order function, thus showing that the finite order
restriction in Corollaries \ref{C:discrete-quotient} and
\ref{T:discrete-quotient-2} cannot be removed. The following
example shows that the exponent ``$\sigma-1+\varepsilon$" that
appears in (\ref{E:discrete-proximity}) cannot be replaced by
``$\sigma-1$".

\begin{example} Since the order of $\Gamma (z)$ is 1, and that
\[
    m\bigg(r,\,\frac{\Gamma(z+1)}{\Gamma(z)}\bigg)=\log r,
\]
we thus see immediately that we cannot drop the $\varepsilon>0$
from (\ref{E:discrete-proximity}). More generally, let $\sigma>0$,
then according to \cite[\S 6, Theorem 5]{Whit} for any given
meromorphic function $\Psi(z)$ of order $\sigma$ there is a
meromorphic solution $F(z)$ to the equation (\ref{E:1st-order})
with  $\sigma(F)\le \sigma +1$. In particular, we may choose
$\Psi(z)$ so that
\begin{equation}
\label{E:Psi}
    m(r,\, \Psi)=T(r,\, \Psi)\sim r^{\sigma}\log
    r.
\end{equation}
Thus we see that
\[
%   \begin{split}
    m\bigg(r,\,\frac{F(z+1)}{F(z)}\bigg) =m(r,\,\Psi) \sim r^{\sigma}\log
    r
    > Cr^{\sigma(F)-1},
%    &=O\big(r^{\sigma(F)-1+\varepsilon}\big).
%    \end{split}
\]
for each positive constant $C$ when we choose $r$ to be
sufficiently large, since $\sigma(F)\le \sigma +1$ holds. We
conclude that we cannot drop the $\varepsilon>0$ from
(\ref{E:discrete-proximity}), and hence
(\ref{E:discrete-proximity}) and (\ref{E:discrete-quotient-2}) are
the best possible in this sense.
\end{example}

\begin{remark} Let $f(z)$ be meromorphic of finite order $\sigma$. Let  $\varepsilon>0$ be given. Then
$T(r,\,f)<O(r^{\sigma+\varepsilon})$.  If we  choose $R=3r,\,
R^\prime=4r$ in (\ref{E:fund-est-I}), then we
obtain
\begin{equation}
\label{E:discrete-quotient-uniform}
    m\left(r,\,\frac{f(z+\eta)}{f(z)}\right)+m\left(r,\,
    \frac{f(z)}{f(z+\eta)}\right)=O(r^{\sigma+\varepsilon}),
\end{equation}
holds \textit{uniformly} for $|\eta|<r$.
\end{remark}

\begin{remark} Let $f(z)$ be a meromorphic function. We choose
$R=3r,\,R^\prime=4r$ in (\ref{E:fund-est-I}). Then, we have
\begin{equation}
\label{E:discrete-est-II}
%\begin{split}
    m\big(r,\, f(z+\eta)\big)\le m(r,\,
    f)+m\left(r,\,\frac{f(z+\eta)}{f(z)}\right)
    = O(T(4r,\, f)).
%\end{split}
\end{equation}
for $|\eta|<r$ to hold \textit{uniformly}.
\end{remark}
It is instructive to compare the (\ref{E:discrete-est-II}) and the
stronger estimate
\begin{equation*}
m\big(r,\,{f(z+\eta)}\big)\le 10\, T(4r,\, f)
\end{equation*}
holds \textit{uniformly} for $|\eta|<r$ for all $r$ sufficiently
large. It is obtained by one of the authors and Ruijsenaars in
\cite[Lem. 3.2]{CR}, by computing directly on the Poisson-Jensen
formula \cite{Hay64}.
\medskip

\section{Some preliminary results}

\begin{lemma}
\label{L:lemma-1}
    Let $\alpha$ be a given constant with $0<\alpha\le 1$. Then
    there exists a constant $C_\alpha>0$ depending only on
    $\alpha$ such that
\begin{equation}
\label{E:lemma-1}
    \log(1+x)\le C_\alpha x^\alpha,
\end{equation}
holds for $x\ge 0$. In particular, $C_1=1$.
\end{lemma}

\begin{proof} The case when $\alpha=1$ is well-known. For $\alpha$
with $0<\alpha<1$, we define the function
\begin{equation*}
     g_{\alpha}(x)=\frac{\log(1+x)}{x^{\alpha}}.
\end{equation*}
It is clear that $g_{\alpha}(x)$ is continous on $(0,\ +\infty)$.
Since
\begin{equation*}
     \lim_{x\rightarrow 0}g_{\alpha}(x)=0,\
    \lim_{x\rightarrow +\infty}g_{\alpha}(x)=0
\end{equation*}
hold. We deduce that $g_{\alpha}$ is bounded on $(0,\ +\infty)$.
So there exists a constant $C_\alpha$,
\begin{equation}
\label{E:lemma-2-constant}
 C_{\alpha}=\max_{0<x<+\infty}g_{\alpha}(x)
\end{equation}
depending only on $\alpha$ such that (\ref{E:lemma-1}) holds.
\end{proof}

\begin{lemma}
\label{L:lemma-2}
    Let $\alpha,\, 0<\alpha\le 1$ be given and $C_\alpha$ as given
    in (\ref{E:lemma-2-constant}). Then for any two complex
    numbers $z_1$ and $z_2$, we have the inequality
\begin{equation}
\label{E:lemma-2}
    \left|\log\left|\frac{z_1}{z_2}\right|\right|\le C_\alpha
    \left(\left|\frac{z_1-z_2}{z_2}\right|^\alpha+
    \left|\frac{z_2-z_1}{z_1}\right|^\alpha\right).
\end{equation}
\end{lemma}

\begin{proof} We deduce from (\ref{E:lemma-1}) that
\begin{equation}
    \log\left|\frac{z_1}{z_2}\right|=
    \log\left|1+\frac{z_1-z_2}{z_2}\right|\leq
    \log\left(1+\left|\frac{z_1-z_2}{z_2}\right|\right)
    \leq  C_{\alpha}\left|\frac{z_1-z_2}{z_2}\right|^{\alpha},
\end{equation}
 and similarly
\begin{equation}
    -\log\left|\frac{z_1}{z_2}\right|=\log\left|\frac{z_2}{z_1}\right|
    =\log\left|1+\frac{z_2-z_1}{z_1}\right|
    \leq\log\left(1+\left|\frac{z_2-z_1}{z_1}\right|\right)
    \leq C_{\alpha}\left|\frac{z_2-z_1}{z_1}\right|^{\alpha}.
\end{equation}
Combining the above two inequalities, we deduce
\begin{equation*}
    \left|\log\left|\frac{z_1}{z_2}\right|\right|
    =\max\left\{\log\left|\frac{z_1}{z_2}\right|,
    -\log\left|\frac{z_1}{z_2}\right|\right\}
    \leq C_{\alpha}\left(\left|\frac{z_1-z_2}{z_2}\right|^{\alpha}
    +\left|\frac{z_2-z_1}{z_1}\right|^{\alpha}\right)
\end{equation*}
as required.
\end{proof}

We need the following result which can be found in \cite[p.
62]{HX} and \cite[p. 66]{JV}.

\begin{lemma}
\label{L:lemma-3}
 Let $\alpha,\, 0<\alpha<1$ be given, then for
every given complex number $w$, we have
\begin{equation}
\label{E:lemma-3}
    \frac{1}{2\pi}\int_0^{2\pi}\frac{1}{|re^{i\theta}-w|^{\alpha}}\,
    d\theta\leq \frac{1}{(1-\alpha)r^{\alpha}}.
\end{equation}
\end{lemma}
\bigskip

\begin{lemma}[\cite{Car}; see also \cite{Lev}]
\label{L:Cartan} Let $z_1, z_2, \cdots, z_p$ be any finite
collection of complex numbers, and let $B>0$ be any given positive
number. Then there exists a finite collection of closed disks
$D_1, D_2, \cdots, D_q$ with corresponding radii $r_1, r_2,
\cdots, r_q$ that satisfy
\begin{equation*}
    r_1+r_2+\cdots+r_q=2B,
\end{equation*}
such that if $z\notin D_j$ for $j=1,2,\cdots,q$, then there is a
permutation of the points $z_1,z_2,\cdots,z_p$, say,
$\hat{z}_1,\hat{z}_2,\cdots,\hat{z}_p$, that satisfies
\begin{equation*}
    |z-\hat{z}_l|>B\frac{l}{p}, \quad\quad l=1, 2, \cdots, p,
\end{equation*}
where the permutation may depend on $z$.
\end{lemma}

\begin{lemma}[A. Z. Mohon'ko \cite{Moho}; see also Laine \cite{Laine}]
\label{L:Mohonko} Let $f(z)$ be a meromorphic function. Then for
all irreducible rational functions in $f$,
\begin{equation}
\label{E:Mohonko}
    R(z,\, f)=\frac{P(z,\, f)}{Q(z,\, f)}=\frac{\sum_{i=0}^p
    a_i(z)f^i}{\sum_{j=0}^q b_j(z)f^j},
\end{equation}
such that the meromorphic coefficients $a_i(z),\, b_j(z)$ satisfy
\begin{equation}
\label{E:Mohonko-2}
    \begin{cases}
    T(r, a_i)=S(r,\, f), & i=0,\,1,\,\cdots, p;\\
    T(r, b_j)=S(r,\, f), & j=0,\,1,\,\cdots, q,
    \end{cases}
\end{equation}
then we have
\begin{equation}
\label{E:Mohonko-3}
    T\big(r,\, R(z,f)\big)=\max\{p,\,q\}\cdot T(r, f)+S(r\, f).
\end{equation}
\end{lemma}
%\medskip
%\vspace{-.2cm}
\section{Proof of Theorem \ref{T:char}}
\begin{proof} Since $f(z)$ has finite order $\sigma$ so that
$\lambda(1/f)\le \sigma<+\infty$. We deduce from Theorem
\ref{T:counting} that
\[
    %\label{E:counting-est}
    N\big(r,f(z+\eta)\big)=N(r,\, f)+O(r^{\lambda-1+\varepsilon})+O(\log r)
\]
holds for the function $f(z)$. This relation and
(\ref{E:discrete-proximity})  together yield
\begin{align*}
    T\big(r,\,f(z+\eta)\big) &= m\big(r, f(z+\eta)\big)+N\big(r,\, f(z+\eta)\big)\\
                     &\le m(r,\, f)+m\Big(r,\,\frac{f(z+\eta)}{f(z)}\Big)
                     + N(r,\,f)+O(r^{\sigma-1+\varepsilon})+O(\log r)\\
                     &=T(r,\,f)+O(r^{\sigma-1+\varepsilon})+O(\log
                     r).
\end{align*}
Similarly, we deduce
\begin{align*}
    T(r,\, f) &= m(r,\,f)+N(r,\, f)\\
                 &\le m\big(r,\, f(z+\eta)\big)+N(r,\, f)+m\Big(r,\,
                 \frac{f(z)}{f(z+\eta)}\Big)+O(\log r)\\
                 &= T\big(r,\, f(z+\eta)\big)+O(r^{\sigma-1+\varepsilon})+O(\log
                     r).
\end{align*}
This completes the proof.
\end{proof}
%\vspace{-.3cm}
\section{Proof of Theorem \ref{T:counting}}
\begin{proof}
 Let $(b_{\mu})_{\mu\in N}$ be the sequence of poles of $f$, with due count of
 multiplicity. Then $(b_{\mu}-\eta)_{\mu\in
N}$ is the sequence of poles of $f(z+\eta)$. Thus by appealing to
the definition of $N(r,\, f)$, we deduce
\begin{equation}
\label{E:proof-thm-count-1}
\begin{split}
    \big|N(r,&\,f(z +\eta))-N(r,f)\big|\\
    &=\bigg|\sum_{0<|b_{\mu}-\eta|<r}\log\frac{r}{|b_{\mu}-\eta|}+n(0,
    f(z+\eta))\log r-\sum_{0<|b_{\mu}|<r}\log\frac{r}{|b_{\mu}|}-n(0, f)\log r\bigg|\\
    &\leq\Bigg|\sum_{
    \substack { 0<|b_{\mu}-\eta|<r,\\
    0<|b_{\mu}|<r }}
    \left(\log\frac{r}{|b_{\mu}-\eta|}-\log\frac{r}{|b_{\mu}|}\right)\Bigg|
    +\Bigg(\sum_{\substack{
    0<|b_{\mu}-\eta|<r,\\ |b_{\mu}|\geq r\ \text{or}\ b_{\mu}=0
    }} \log\frac{r}{|b_{\mu}-\eta|}\Bigg)\\
    &\quad\quad +\Bigg(\sum_{
    \substack{
    0<|b_{\mu}|<r,\\
     |b_{\mu}-\eta|\geq r\ \text{or}\ b_{\mu}-\eta=0
    }} \log\frac{r}{|b_{\mu}|}\Bigg)+O(\log r)\\
    &\leq \Bigg(\sum_{
    \substack{
    0<|b_{\mu}-\eta|<r,\\ 0<|b_{\mu}|<r}}
    \bigg|\log\bigg|\frac{b_{\mu}}{b_{\mu}-\eta}\bigg|\bigg|\Bigg)
    +\bigg(\sum_{
    \substack{
    0<|b_{\mu}-\eta|<r,\\ |b_{\mu}|\geq r}}
    \log\frac{r}{|b_{\mu}-\eta|}\Bigg)
    +\Bigg(\sum_{
    \substack{
    0<|b_{\mu}|<r,\\ |b_{\mu}-\eta|\geq r}}
    \log\frac{r}{|b_{\mu}|}\Bigg)\\
    &\quad\quad+O(\log r).
    \end{split}
    \end{equation}
Applying Lemma \ref{L:lemma-2} with $\alpha=1$ to an individual
term in the first summand of the last inequality
(\ref{E:proof-thm-count-1}), we deduce
\begin{equation}
\label{E:proof-thm-count-2}
\begin{split}
    \bigg|\log\bigg|\frac{b_{\mu}}{b_{\mu}-\eta}\bigg|\bigg| &\le
    \bigg|\frac{b_\mu-(b_\mu-\eta)}{b_\mu-\eta}\bigg|+
    \bigg|\frac{(b_\mu-\eta)-b_\mu}{b_\mu}\bigg|\\
    &=\bigg|\frac{\eta}{b_\mu-\eta}\bigg|+\bigg|\frac{\eta}{b_\mu}\bigg|.
\end{split}
\end{equation}
We now consider the second summand in the last line of
(\ref{E:proof-thm-count-1}). More specifically, we apply the Lemma
\ref{L:lemma-1} and inequalities $0<|b_\mu-\eta|< r,\, |b_\mu|\ge
r$ that restrict the summation to obtain the inequalities
\begin{equation}
\label{E:proof-thm-count-3}
\begin{split}
    \log\frac{r}{|b_\mu-\eta|}&=\log\bigg(\frac{r-|b_\mu-\eta|}{|b_\mu-\eta|}+1\bigg)
    \le\frac{r-|b_\mu-\eta|}{|b_\mu-\eta|}\\
    &\le\frac{|\eta|+r-|b_\mu|}{|b_\mu-\eta|}\le\frac{|\eta|}{|b_\mu-\eta|}.
\end{split}
\end{equation}
Let us consider the third summand in the last line of
(\ref{E:proof-thm-count-1}). We similarly consider the
inequalities $|b_\mu-\eta|\ge r,\,|b_\mu|< r$ that restrict the
summation. This gives $|\eta|\ge r-|b_\mu|>0$. We conclude from
Lemma \ref{L:lemma-1} the inequalities
\begin{equation}
\label{E:proof-thm-count-4}
%\begin{split}
    \log\frac{r}{|b_\mu|} =
    \log\bigg(\frac{r-|b_\mu|}{|b_\mu|}+1\bigg)\le \frac{r-|b_\mu|}{|b_\mu|}
    \le\bigg|\frac{\eta}{b_\mu}\bigg|.
%\end{split}
\end{equation}
Combining the above inequalities (\ref{E:proof-thm-count-1}),
(\ref{E:proof-thm-count-2}), (\ref{E:proof-thm-count-3}) and
(\ref{E:proof-thm-count-4}), we deduce
\begin{equation}
\label{E:proof-thm-count-5}
\begin{split}
\big|N(r,&\,f(z +\eta))-N(r,f)\big|\\
    &\leq |\eta|\Bigg\{\sum_{
    \substack{
    0<|b_{\mu}-\eta|<r,\\ 0<|b_{\mu}|<r}}
    \left(\frac{1}{|b_{\mu}|}+\frac{1}{|b_{\mu}-\eta|}\right)
    +\sum_{
    \substack{
    0<|b_{\mu}-\eta|<r,\\ |b_{\mu}|\geq r}}
    \frac{1}{|b_{\mu}-\eta|}
    +\sum_{
    \substack{
    0<|b_{\mu}|<r,\\ |b_{\mu}-\eta|\geq r}}
    \frac{1}{|b_{\mu}|}\Bigg\}+O(\log r)\\
    &= |\eta|\left\{\sum_{0<|b_{\mu}-\eta|<r}\frac{1}{|b_{\mu}-\eta|}
    +\sum_{0<|b_{\mu}|<r}\frac{1}{|b_{\mu}|}\right\}+O(\log r).
\end{split}
\end{equation}
We turn to estimate the first summand in the last line of
(\ref{E:proof-thm-count-5}). In particular, we divide the
summation range $0<|b_{\mu}-\eta|<r$ into two ranges, namely the
$0<|b_{\mu}-\eta|\le |\eta|$ and $|\eta|<|b_{\mu}-\eta|<r$. We
notice that when $|b_{\mu}-\eta|>|\eta|$, then
\begin{equation}
\label{E:proof-thm-count-6}
    \frac{1}{|b_{\mu}-\eta|}=\frac{1}{|b_{\mu}|}\cdot\bigg|1+\frac{\eta}{b_{\mu}-\eta}\bigg|
    \leq \frac{1}{|b_{\mu}|}\cdot\bigg(1+\bigg|\frac{\eta}{b_{\mu}-\eta}\bigg|\bigg)<
    \frac{2}{|b_{\mu}|}.
\end{equation}
Thus when $r>|\eta|$ the first summand on the last line of
(\ref{E:proof-thm-count-5}) becomes
\begin{equation}
\label{E:proof-thm-count-7}
\begin{split}
    \sum_{0<|b_{\mu}-\eta|<r}\frac{1}{|b_{\mu}-\eta|}=&
    \sum_{0<|b_{\mu}-\eta|\leq|\eta|}\frac{1}{|b_{\mu}-\eta|}
    +\sum_{|\eta|<|b_{\mu}-\eta|<r}\frac{1}{|b_{\mu}-\eta|}\\
    \leq &
    2\cdot\bigg(\sum_{|\eta|<|b_{\mu}-\eta|<r}\frac{1}{|b_{\mu}|}\bigg)+O(1)\\
    \leq & 2\cdot\bigg(\sum_{0<|b_{\mu}|<r+|\eta|}\frac{1}{|b_{\mu}|}\bigg)+O(1).
\end{split}
\end{equation}
Combining (\ref{E:proof-thm-count-5})and
(\ref{E:proof-thm-count-7}), we get
\begin{equation}
\label{E:proof-thm-count-8}
    \Big|N(r,f(z+\eta))-N(r, f)\Big|\leq 3|\eta|
    \bigg(\sum_{0<|b_{\mu}|<r+|\eta|}\frac{1}{|b_{\mu}|}\bigg)+O(\log r).
\end{equation}

We distinguish two cases:
\begin{enumerate}
\item{Case 1:} $\lambda\geq 1$. By the H\"{o}lder inequality, we
have for any $\varepsilon>0$,
\begin{equation}
\label{E:proof-thm-count-9}
\begin{split}
    \sum_{0<|b_{\mu}|<r+|\eta|}\frac{1}{|b_{\mu}|} &\leq
    \bigg(\sum_{0<|b_{\mu}|<r+|\eta|}\frac{1}{|b_{\mu}|^
    {\lambda+\varepsilon}}\bigg)^{\frac{1}{\lambda+\varepsilon}}
    \cdot\bigg(\sum_{0<|b_{\mu}|<r+|\eta|}1^{\frac{\lambda+\varepsilon}{\lambda+\varepsilon-1}}
    \bigg)^{\frac{\lambda+\varepsilon-1}{\lambda+\varepsilon}}\\
    &\le O(1)\cdot n(r+|\eta|,f)^{\frac{\lambda+\varepsilon-1}{\lambda+\varepsilon}}.
\end{split}
\end{equation}
But
\begin{equation}
\label{E:proof-thm-count-10}
    n(r+|\eta|,f)=O((r+|\eta|)^{\lambda+\varepsilon})=O(r^{\lambda+\varepsilon}).
\end{equation}
Therefore, inequalities (\ref{E:proof-thm-count-9}) and
(\ref{E:proof-thm-count-10}) give
\begin{equation}
\label{E:proof-thm-count-11}
\sum_{0<|b_{\mu}|<r+|\eta|}\frac{1}{|b_{\mu}|}=O(r^{\lambda-1+\varepsilon}).
\end{equation}

\item{Case 2:} $\lambda<1$. We have, by the definition of exponent
of convergence,
\begin{equation}
\label{E:proof-thm-count-12}
\sum_{0<|b_{\mu}|<r+|\eta|}\frac{1}{|b_{\mu}|}=O(1).
\end{equation}
\end{enumerate}
We finally obtain from (\ref{E:proof-thm-count-8}),
(\ref{E:proof-thm-count-11}) and (\ref{E:proof-thm-count-12}) the
desired result
\begin{equation*}
    |N\big(r,f(z+\eta)\big)-N(r, f)|=O(r^{\lambda-1+\varepsilon})+O(\log r).
\end{equation*}
\end{proof}

\section{Proof of Theorem \ref{T:counting-example}}
\begin{proof}
Let $\alpha,\, 0<\alpha\le 1$, and let a sequence of numbers
located at positive integers $k,\, k=2,\, 3,\,4,\,  \cdots$, each
with multiplicity $\gamma_k$. Then according to Weierstrass'
theorem \cite[\S8.1]{Titch}, there is an entire function $g(z)$
that has zeros precisely at the sequence defined above. We now
take $f(z)=1/g(z)$ to be the meromorphic function that we consider
below. We then write
\begin{equation}
\label{E:proof-count-2-2}
    N(r,\, f)=\sum_{2\le k<r}\gamma_k \log\frac{r}{k}.
\end{equation}
Since the poles of $f(z+1)$ are those of $f(z)$ but shifted to the
left by one unit, so let us write
\[
    N\big(r,\, f(z+1)\big)=\sum_{1\le k<r}\beta_k \log\frac{r}{k},
\]
where $\beta_k=\gamma_{k+1}$ for $k=1,\,2,\, 3,\, 4,\,\cdots$. We
deduce
\begin{equation}
\label{E:proof-count-2-3}
\begin{split}
    N\big(r,f(z+1)\big)-N(r, f) &=\sum_{1\le k<r}\beta_{k} \log\frac{r}{k}
    -\sum_{2\le k<r}\gamma_k \log\frac{r}{k}\\
    &=\sum_{1\le k<r}\gamma_{k+1} \log\frac{r}{k}
    -\sum_{2\le k<r}\gamma_k \log\frac{r}{k}\\
    &=\gamma_2 \log r +\sum_{2\le k<r} \big(\gamma_{k+1}-\gamma_k\big)
    \log\frac{r}{k}.
\end{split}
\end{equation}
We now choose
\begin{equation}
\label{E:proof-count-2-1}
    \gamma_{k+1}=2\gamma_k,\quad k=2,\, 3,\, 4,\, \cdots
\end{equation}
then
\begin{equation}
\label{E:proof-count-2-5}
\begin{split}
    \frac{N\big(r,\, f(z+\eta)\big)-N\big(r,\, f(z)\big)}{N\big(r,\, f(z)\big)}
    &\ge
    \frac{\gamma_2 \log r +\sum_{2\le k<r} \big(\gamma_{k+1}-\gamma_k\big)
    \log\frac{r}{k}}{\sum_{2\le k<r}\gamma_k \log\frac{r}{k}}\\
    &= 1+\frac{\gamma_2 \log r}{\sum_{2\le k<r}\gamma_k \log\frac{r}{k}}\ge 1
%    &=1+\frac{\gamma_2 \log r}{N(r,\,f)}= 1+o(1),
\end{split}
\end{equation}
for all $r\ge 3$. On the other hand, it is easy to see from the
meromorphic function $g(z)$ constructed above that it has an
infinite order of growth.
\end{proof}

\section{Proof of Theorem \ref{T:discrete-quotient}}
\begin{proof}
Let $z=re^{i\theta}$ such that $|z|<R-|\eta|$. The Poisson-Jensen
formula yields
\begin{equation}
\label{E:proof-thm-prox-1}
\begin{split}
    \log |f(z)|=&\frac{1}{2\pi}\int_0^{2\pi}\log |f(Re^{i\phi})|
    \Re\left(\frac{Re^{i\phi}+z}{Re^{i\phi}-z}\right)\,d\phi\\
                &-\sum_{|a_{\nu}|<R}\log\left|\frac{R^2-\bar{a}_{\nu}z}{R(z-a_{\nu})}\right|+
                \sum_{|b_{\mu}|<R}\log\left|\frac{R^2-\bar{b}_{\mu}z}{R(z-b_{\mu})}\right|,
\end{split}
\end{equation}
where $(a_{\nu})_{\nu\in N}$ and $(b_{\mu})_{\mu\in N}$, denote
respectively, and with due count of multiplicity, the zeros and
poles of $f$ in $\{|z|<R\}$. Since $|z+\eta|<R$, so the
(\ref{E:proof-thm-prox-1}) also yields
\begin{equation}
\label{E:proof-thm-prox-2}
\begin{split}
    \log |f(z+\eta)|=&\frac{1}{2\pi}\int_0^{2\pi}\log
    |f(Re^{i\phi})|\,
    \Re\left(\frac{Re^{i\phi}+z+\eta}{Re^{i\phi}-z-\eta}\right)\, d\phi\\
    &-\sum_{|a_{\nu}|<R}\log\left|\frac{R^2-\bar{a}_{\nu}(z+\eta)}{R(z+\eta-a_{\nu})}\right|+
    \sum_{|b_{\mu}|<R}\log\left|\frac{R^2-\bar{b}_{\mu}(z+\eta)}{R(z+\eta-b_{\mu})}\right|.
\end{split}
\end{equation}

Subtracting (\ref{E:proof-thm-prox-1}) from
(\ref{E:proof-thm-prox-2}) yields
\begin{equation}
\label{E:proof-thm-prox-3}
\begin{split}
    \log \left|\frac{f(z+\eta)}{f(z)}\right|= &\frac{1}{2\pi}\int_0^{2\pi}
    \log\left|f(Re^{i\phi})\right|\,\Re\left(\frac{2\eta Re^{i\phi}}
    {(Re^{i\phi}-z-\eta)(Re^{i\phi}-z)}\right)\,d\phi\\
    &-\sum_{|a_{\nu}|<R}\log\left|\frac{R^2-\bar{a}_{\nu}(z+\eta)}{R^2-\bar{a}_{\nu}z}\right|+
    \sum_{|b_{\mu}|<R}\log\left|\frac{R^2-\bar{b}_{\mu}(z+\eta)}{R^2-\bar{b}_{\mu}z}\right|\\
    &+\sum_{|a_{\nu}|<R}\log\left|\frac{R(z+\eta-a_{\nu})}{R(z-a_{\nu})}\right|-
    \sum_{|b_{\mu}|<R}\log\left|\frac{R(z+\eta-b_{\mu})}{R(z-b_{\mu})}\right|.
\end{split}
\end{equation}
We deduce from (\ref{E:proof-thm-prox-3}) that
\begin{equation}
\label{E:proof-thm-prox-4}
\begin{split}
    \left|\log\left|\frac{f(z+\eta)}{f(z)}\right|\right|\leq
    &\left(\frac{2|\eta| R}{(R-|z|-|\eta|)(R-|z|)}\right)\cdot\frac{1}{2\pi}\int_0^{2\pi}\big|
    \log |f(Re^{i\phi})|\big|\,d\phi\\
    &+\sum_{|a_{\nu}|<R}\left|\log\left|\frac{R^2-\bar{a}_{\nu}(z+\eta)}{R^2-\bar{a}_{\nu}z}
    \right|\right|+\sum_{|b_{\mu}|<R}\left|\log\left|\frac{R^2-\bar{b}_{\mu}(z+\eta)}
    {R^2-\bar{b}_{\mu}z}\right|\right|\\
    &+\sum_{|a_{\nu}|<R}\left|\log\left|\frac{z+\eta-a_{\nu}}{z-a_{\nu}}\right|\right|+
    \sum_{|b_{\mu}|<R}\left|\log\left|\frac{z+\eta-b_{\mu}}{z-b_{\mu}}\right|\right|\\
    \leq &\frac{2|\eta| R}{(R-|z|-|\eta|)^2}\big(m(R,f)+m(R,{1}/{f})\big)\\
    &+\sum_{|a_{\nu}|<R}\left|\log\left|\frac{R^2-\bar{a}_{\nu}(z+\eta)}
    {R^2-\bar{a}_{\nu}z}\right|\right|+\sum_{|b_{\mu}|<R}
    \left|\log\left|\frac{R^2-\bar{b}_{\mu}(z+\eta)}{R^2-\bar{b}_{\mu}z}\right|\right|\\
    &+\sum_{|a_{\nu}|<R}\left|\log\left|\frac{z+\eta-a_{\nu}}{z-a_{\nu}}\right|\right|+
    \sum_{|b_{\mu}|<R}\left|\log\left|\frac{z+\eta-b_{\mu}}{z-b_{\mu}}\right|\right|.
\end{split}
\end{equation}

We apply (\ref{E:lemma-2}) with $\alpha=1$ and Lemma
\ref{L:lemma-2} to the second and third summands in
(\ref{E:proof-thm-prox-4}). This yields, for $|a_{\nu}|<R$,
\begin{equation}
\label{E:proof-thm-prox-5}
\begin{split}
    \left|\log\left|\frac{R^2-\bar{a}_{\nu}(z+\eta)}{R^2-\bar{a}_{\nu}z}\right|\right|
    \leq & \left|\frac{\bar{a}_{\nu}\eta}{R^2-\bar{a}_{\nu}z}\right|+
    \left|\frac{\bar{a}_{\nu}\eta}{R^2-\bar{a}_{\nu}(z+\eta)}\right|\\
    \leq & \frac{|\eta|}{R-|z|}+\frac{|\eta|}{R-|z|-|\eta|}\leq
    \frac{2|\eta|}{R-|z|-|\eta|}.
\end{split}
\end{equation}
Similarly, we have for $|b_{\mu}|<R$,
\begin{equation}
\label{E:proof-thm-prox-6}
\left|\log\left|\frac{R^2-\bar{b}_{\mu}(z+\eta)}{R^2-\bar{b}_{\mu}z}\right|\right|
\leq \frac{2|\eta|}{R-|z|-|\eta|}.
\end{equation}

We then choose $0<\alpha<1$ in (\ref{E:lemma-2}) and this yields
\begin{equation}
\label{E:proof-thm-prox-7}
    \left|\log\left|\frac{z+\eta-a_{\nu}}{z-a_{\nu}}\right|\right|\leq
    C_{\alpha}|\eta|^{\alpha}\left(\frac{1}{|z-a_{\nu}|^{\alpha}}
    +\frac{1}{|z+\eta-a_{\nu}|^{\alpha}}\right),
\end{equation}
and
\begin{equation}
\label{E:proof-thm-prox-8}
    \left|\log\left|\frac{z+\eta-b_{\mu}}{z-b_{\mu}}\right|\right|\leq
    C_{\alpha}|\eta|^{\alpha}\left(\frac{1}{|z-b_{\mu}|^{\alpha}}
    +\frac{1}{|z+\eta-b_{\mu}|^{\alpha}}\right).
\end{equation}

Combining the inequalities
(\ref{E:proof-thm-prox-4})--(\ref{E:proof-thm-prox-8}), we get
\begin{equation}
\label{E:proof-thm-prox-9}
\begin{split}
    \left|\log\left|\frac{f(z+\eta)}{f(z)}\right|\right|\leq &\frac{2|\eta|R}
    {(R-|z|-|\eta|)^2}\bigg(m\big(R,\, {f}\big)+m\big(R,\,\frac{1}{f}\big)\bigg)\\
    &+\frac{2|\eta|}{R-|z|-|\eta|}\bigg(n\big(R,\, {f}\big)+n\big(R,\,\frac{1}{f}\big)\bigg)\\
    & +C_{\alpha}|\eta|^{\alpha}\sum_{|a_{\nu}|<R}\left(\frac{1}{|z-a_{\nu}|^{\alpha}}
    +\frac{1}{|z+\eta-a_{\nu}|^{\alpha}}\right)\\
    &+C_{\alpha}|\eta|^{\alpha}\sum_{|b_{\mu}|<R}
    \left(\frac{1}{|z-b_{\mu}|^{\alpha}}+\frac{1}{|z+\eta-b_{\mu}|^{\alpha}}\right).
\end{split}
\end{equation}
Integrating (\ref{E:proof-thm-prox-9}) on $|z|=r$, and applying
Lemma \ref{L:lemma-3} gives
\begin{equation}
\label{E:proof-thm-prox-10}
\begin{split}
    m\left(r,\,\frac{f(z+\eta)}{f(z)}\right)&+
    m\left(r,\,\frac{f(z)}{f(z+\eta)}\right)\\
    &\leq  \frac{2|\eta|R}{(R-r-|\eta|)^2}\bigg(m\big(R,\, {f}\big)+m\big(R,\,\frac{1}{f}\big)\bigg)
     +\frac{2|\eta|}{R-r-|\eta|}\big(n(R,\,f)+n(R,\,{1}/{f})\big)\\
    &\quad +C_{\alpha}|\eta|^{\alpha}\sum_{|a_{\nu}|<R}\left(\frac{1}{2\pi}\int_0^{2\pi}
    \frac{1}{|re^{i\theta}-a_{\nu}|^{\alpha}}\,d\theta+\frac{1}{2\pi}\int_0^{2\pi}
    \frac{1}{|re^{i\theta}+\eta-a_{\nu}|^{\alpha}}\,d\theta\right)\\
    &\quad +C_{\alpha}|\eta|^{\alpha}\sum_{|b_{\mu}|<R}\left(\frac{1}{2\pi}\int_0^{2\pi}
    \frac{1}{|re^{i\theta}-b_{\mu}|^{\alpha}}\,d\theta+\frac{1}{2\pi}\int_0^{2\pi}
    \frac{1}{|re^{i\theta}+\eta-b_{\mu}|^{\alpha}}\,d\theta\right)\\
    &\leq \frac{2|\eta| R}{(R-r-|\eta|)^2}\bigg(m\big(R,\, {f}\big)+m\big(R,\,\frac{1}{f}\big)\bigg)\\
    &\quad +\left(\frac{2|\eta|}{R-r-|\eta|}+\frac{2C_{\alpha}|\eta|^{\alpha}}
    {(1-\alpha)r^{\alpha}}\right)\bigg(n\big(R,\, {f}\big)+n\big(R,\,\frac{1}{f}\big)\bigg).
\end{split}
\end{equation}
Since $R^{\prime}>R>1$, we deduce
\begin{equation*}
\begin{split}
    N(R^{\prime},\, f) &\geq \int_R^{R^{\prime}}\frac{n(t,\,f)-n(0,\, f)}{t}\,dt
    +n(0,\,f)\log R^{\prime}\\
    &\geq n(R,\,f)\int_R^{R^{\prime}}\frac{dt}{t}-n(0,\,f)\int_R^{R^{\prime}}\frac{dt}{t}
    +n(0,f)\log R^{\prime}\\
    &\geq n(R,\,f)\frac{R^{\prime}-R}{R^{\prime}}.
\end{split}
\end{equation*}
Hence
\begin{equation}
\label{E:proof-thm-prox-11}
    n(R,\,f)\leq\frac{R^{\prime}}{R^{\prime}-R}\, N(R^{\prime},\,f).
\end{equation}
Similarly, we have
\begin{equation}
\label{E:proof-thm-prox-12}
    n\big(R,\,\frac{1}{f}\big)\leq\frac{R^{\prime}}{R^{\prime}-R}\, N\big(R^{\prime},\,\frac{1}{f}\big).
\end{equation}
We deduce the (\ref{E:fund-est-I}) after combining
(\ref{E:proof-thm-prox-10}), (\ref{E:proof-thm-prox-11}) and
(\ref{E:proof-thm-prox-12}).
\end{proof}

\section{Pointwise estimates}

It is well-known that pointwise logarithmic derivative estimates
of finite order meromorphic functions play an important role in
complex differential equations (See e.g. \cite{Gund-2}). In
particular, the following estimate of Gundersen \cite[Cor.
2]{Gund} gives a sharp upper bound of logarithm derivatives.

\begin{theorem}
\label{T:Gund-1} Let $f(z)$ be a meromorphic function, and let
$k\ge 1$ be an integer, $\alpha>1$, and $\varepsilon>0$ be given
real constants, then there exists a set $E\subset (1,\infty)$ of
finite logarithmic measure,
\begin{enumerate}
\item [(a)] and a constant $A>0$ depending only on $\alpha$, such
that for all $|z|\not\in E\cup[0,\,1]$, we have
\begin{equation}
\label{E:thm-Gund-a}
    \bigg|\frac{f^{\prime}(z)}{f(z)}\bigg|\le A\bigg(\frac{T(\alpha
    r,f)}{r}+\frac{n(\alpha r)}{r}\log^\alpha r\log^+n(\alpha
    r)\bigg),
\end{equation}
where $n(t)=n(t,f)+n(t,1/f)$;

\item [(b)] and if in addition that $f(z)$ has finite order
$\sigma$, and such that for all $|z|\not\in E\cup[0,\,1]$, we have
\begin{equation}
\label{E:thm-Gund-b}
 \bigg|\frac{f^{\prime}(z)}{f(z)}\bigg|\le
    |z|^{\sigma-1+\varepsilon}.
\end{equation}
\end{enumerate}
\end{theorem}

We first give pointwise estimates for our difference quotient
which are counterparts to Gundersen's logarithmic derivative
estimates.

\begin{theorem}
\label{T:thm-ptwse-2} Let $f(z)$ be a meromorphic function, $\eta$
a non-zero complex number, and let $\gamma>1$, and $\varepsilon>0$
be given real constants, then there exists a subset $E\subset
(1,\infty)$ of finite logarithmic measure,
\begin{enumerate}
\item[(a)] and a constant $A$ depending only on $\gamma$ and
$\eta$, such that for all $|z|\not\in E\cup[0,\,1]$, we have
\begin{equation}
\label{E:thm-ptwse-2-a}
    \left|\log\bigg|\frac{f(z+\eta)}{f(z)}\bigg|\right|\le A\left(
    \frac{T(\gamma r, f)}{r}+\frac{n(\gamma r)}{r}\log^\gamma r
    \log^+n(\gamma r)\right);
\end{equation}

\item[(b)] and if in addition that $f(z)$ has finite order
$\sigma$, and such that for all \break
    $|z|=r\not\in E\cup[0,\,1]$, we
have
\begin{equation}
\label{E:thm-ptwse-2-b}
    \exp\big(-r^{\sigma-1+\varepsilon}\big)\leq
    \bigg|\frac{f(z+\eta)}{f(z)}\bigg|\leq
    \exp\big(r^{\sigma-1+\varepsilon}\big).
\end{equation}
\end{enumerate}
\end{theorem}
\bigskip

We remark that the example $f(z)=e^{z^n}$ shows that the
$\varepsilon>0$ in the (\ref{E:thm-ptwse-2-b}) cannot be dropped,
and so the (\ref{E:thm-ptwse-2-b}) is the best possible.

 The forms of logarithmic derivative estimates almost always depend on how
we remove the ``exceptional set" consisting of the zeros and poles
of the function in the complex plane, such as the proof given by
Hille \cite[Thm. 4.5.1]{Hil}. More precise estimates usually
depend on application of the Cartan lemma \cite{Car} (see also
\cite{Lev}) such as \cite[Prop. 5.12]{Laine} and Theorem
\ref{T:Gund-1} above. We shall make use of the same lemma to prove
our Theorem.

\begin{proof} Let $z$ be such that $|z|=r<R-|\eta|$. We deduce
from (\ref{E:proof-thm-prox-3}), (\ref{E:proof-thm-prox-4}) the
inequality (\ref{E:proof-thm-prox-9}).

 Let $\beta>1$ and $R=\beta r+|\eta|$. We choose $r_1$ in
 (\ref{E:proof-thm-prox-9}) so that $|\eta|< \beta(\beta-1)r$ for $r> r_1$. In addition, we
apply Lemma \ref{L:lemma-2} with $\alpha=1$ (note that $C_1=1$),
and this yields
\begin{equation}
\label{E:proof-thm-ptwse-2-2}
\begin{split}
    \bigg|\log\bigg|\frac{f(z+\eta)}{f(z)}\bigg|\bigg| &\leq
    \frac{4|\eta|(\beta r+|\eta|)}{(\beta-1)^2r^2} T\big(\beta
    r+|\eta|, f\big) +|\eta|\cdot\sum_{|c_k|<\beta r+|\eta|}
    \bigg(\frac{1}{|z-c_k|}+\frac{1}{|z+\eta-c_k|}\bigg)\\
    &\le \frac{4|\eta|\beta^2}{(\beta-1)^2}
    \frac{T\big(\beta^2 r,\, f\big)}{r} +|\eta|\cdot\sum_{|c_k|<\beta^2r}
    \bigg(\frac{1}{|z-c_k|}+\frac{1}{|z+\eta-c_k|}\bigg)\\
    &=4|\eta|\left(\frac{\beta}{\beta-1}\right)^2
    \frac{T\big(\beta^2 r,\, f\big)}{r}
    +|\eta|\cdot\sum_{|c_k|<\beta^2r}
    \frac{1}{|z-d_k|},
\end{split}
\end{equation}
where $(c_k)_{k\in N}=(a_{\nu})_{\nu\in N}\cup(b_{\mu})_{\mu\in
N}$ and $(d_k)_{k\in N}=(c_k)_{k\in N}\cup(c_k-\eta)_{k\in N}$,
where the sequence $d_k$ is listed according to multiplicity and
ordered by increasing modulus.

We now let $\gamma=\beta^2$ and apply Lemma \ref{L:Cartan} to the
second summand to (\ref{E:proof-thm-ptwse-2-2}) with
$|d_k|<R=\gamma r$, where $r> \max\{r_1,\, d_1\}$. The argument
then follows the same argument as \cite[(7.6)--(7.9)]{Gund} (with
their $\alpha$ replaced by our $\gamma$) so that we deduce for all
$|z|\not\in E\cup[0,\,1]$, where the $E$ has finite logarithmic
measure,
\begin{equation}
\label{E:proof-thm-ptwse-2-3}
    \sum_{|c_k|<\gamma r}\frac{1}{|z-d_k|}\le
    \gamma^2\frac{n(\gamma^2 r)}{r}\log^\gamma r\log n(\gamma^2 r).
\end{equation}
Combining (\ref{E:proof-thm-ptwse-2-2}) and
(\ref{E:proof-thm-ptwse-2-3}) we obtain
\begin{equation}
\label{E:proof-thm-ptwse-2-4}
\begin{split}
    \bigg|\log\bigg|\frac{f(z+\eta)}{f(z)}\bigg|\bigg| &\le
    4|\eta|\left(\frac{\beta}{\beta-1}\right)^2
    \frac{T\big(\gamma r, f\big)}{r}+|\eta|\cdot
    \bigg(\gamma^2\frac{n(\gamma^2 r)}{r}\log^\gamma r\log n(\gamma^2
    r)\bigg)\\
    &\le
    |\eta|\left[4\left(\frac{\beta}{\beta-1}\right)^2
    \frac{T\big(\gamma^2 r, f\big)}{r} +
    \gamma^2\frac{n(\gamma^2 r)}{r}\log^{\gamma^2} r\log n(\gamma^2
    r)\right]
\end{split}
\end{equation}
which gives (\ref{E:thm-ptwse-2-a}) with $\gamma^2$ replaced by
$\gamma$.

If $f(z)$ has finite order $\sigma$, then given $\varepsilon>0$,
it is now easy to deduce (\ref{E:thm-ptwse-2-b}) holds from the
estimate (\ref{E:thm-ptwse-2-a}).
\end{proof}

We easily obtain the following result.
\begin{corollary}
\label{C:pointwise-1} Let $\eta_1,\, \eta_2$ be two arbitrary
complex numbers, and let $f(z)$ be a meromorphic function of
finite order $\sigma$. Let $\varepsilon>0$ be given, then there
exists a subset $E\subset \textbf{\textbf {R}}$ with finite
logarithmic measure such that for all $r\not\in E\cup [0,\,1]$, we
have
\begin{equation}
\label{E:pointwise-1}
    \exp\big(-r^{\sigma-1+\varepsilon}\big)\leq
    \bigg|\frac{f(z+\eta_1)}{f(z+\eta_2)}\bigg|\leq
    \exp\big(r^{\sigma-1+\varepsilon}\big).
\end{equation}
\end{corollary}
\bigskip

We can replace the linear exceptional set by ``radial exceptional"
set (see also \cite[Lem. 2 and Cor. 4]{Gund} and \cite[Thm. 2 and
Cor. 1]{Gund}).

\begin{theorem}
%\label{T:thm8-6}
Let $f(z)$ be a meromorphic function, $\eta$ a non-zero complex
number, and let $\gamma>1$, and $\varepsilon>0$ be given real
constants, then there exists a set $E\subset [0,2\pi)$ that has
linear measure zero, such that if $z=re^{\psi_0}$ satisfying
$\psi_0\not\in E$, then there is a constant $R_0=R_0(\psi_0)>1$
such that for $z$ satisfying $\arg z=\psi_0$ and $|z|\ge R_0$,
\begin{enumerate}
\item[(a)] we have
\begin{equation*}
%\label{E:thm8-6-a}
    \left|\log\bigg|\frac{f(z+\eta)}{f(z)}\bigg|\right|\le B\left(
    \frac{T(\gamma r, f)}{r}+\frac{n(\gamma r)}{r}\log^\gamma r
    \log^+n(\gamma r)\right),
\end{equation*}
and the constant $B$ depending only on $\gamma$ and $\eta$;

\item[(b)] and if in addition that $f(z)$ has finite order
$\sigma$, we have
\begin{equation*}
%\label{E:thm8-6-b}
    \exp\big(-r^{\sigma-1+\varepsilon}\big)\leq
    \bigg|\frac{f(z+\eta)}{f(z)}\bigg|\leq
    \exp\big(r^{\sigma-1+\varepsilon}\big).
\end{equation*}
\end{enumerate}
\end{theorem}

We shall omit the proof. Similarly we have

\begin{corollary}
%\label{C:pointwise-2}
Let $\eta_1,\, \eta_2$ be two arbitrary complex numbers, and let
$f(z)$ be a meromorphic function of finite order $\sigma$. Let
$\varepsilon>0$ be given, then there exists a subset $E\subset
[0,\,2\pi)$ of linear measure zero  such that if $z=re^{\psi_0}$
satisfying $\psi_0\not\in E$, then there is a constant
$R_0=R_0(\psi_0)>1$ such that for $z$ satisfying $\arg z=\psi_0$
and $|z|\ge R_0$, we have
\begin{equation*}
%\label{E:pointwise-1}
    \exp\big(-r^{\sigma-1+\varepsilon}\big)\leq
    \bigg|\frac{f(z+\eta_1)}{f(z+\eta_2)}\bigg|\leq
    \exp\big(r^{\sigma-1+\varepsilon}\big).
\end{equation*}
\end{corollary}
\medskip

\section{Applications to difference equations}

We first apply the Theorem \ref{T:char} to give a direct proof of
the following theorem.

\begin{theorem}[\cite{AHH}, \cite{HKLRT}]
\label{T:HKLRT}  Let $c_1,\, \cdots, c_n$ be non-zero complex
numbers. If the difference equation
\begin{equation}
\label{E:HKLRT}
    \sum_{i=j}^n y(z+c_j)=R(z,\,y(z))=\frac{a_0(z)+a_1(z)y(z)+\cdots+a_p(z)y(z)^p}
    {b_0(z)+b_1(z)y(z)+\cdots+b_q(z)y(z)^q}
\end{equation}
with polynomial coefficients $a_i,\, b_j$, admits a finite order
meromorphic solution $f(z)$, then we have $\max\{p,\,q\}\le n$.
\end{theorem}

This theorem was first given in \cite[Thm. 3]{AHH} with $n=2$ and
was written in the above generalized form in \cite[Prop.
2.1]{HKLRT}.

\begin{proof}
%This theorem was first given in \cite[Thm. 3]{AHH} with $n=2$ and
%generalized in \cite[Prop. 2.1]{HKLRT} to the above form.
Without loss of generality, we assume $f(z)$ to be a finite order
transcendental meromorphic solution to (\ref{E:HKLRT}). The
estimate on the right side of (\ref{E:HKLRT}) is easily handled by
applying Lemma \ref{L:Mohonko}, as in the proofs in \cite{AHH} and
\cite{HKLRT}, to give (\ref{E:Mohonko-3}).
%$T\big(r,\,R(f)\big)=\max\{p,\,q\}T(r,\,f)$.
Then the (\ref{E:char-est}) of our Theorem \ref{T:char} and
(\ref{E:HKLRT}) yield
\begin{equation}
\label{P:thm9-1}
\begin{split}
    \max\{p,\,q\}T(r,\,f) &=T\big(r,\, R(z,\,f)\big)+S(r,\,f) \\
    &\le T\Big(r,\,\sum_{j=1}^n f(z+c_j)\Big)+S(r,\,f)\\
    &\le nT(r,\,f)+O\big(r^{\sigma-1+\varepsilon}\big)+O(\log r)+
    S(r,\,f)
\end{split}
\end{equation}
since the (\ref{E:char-est}) is independent of $c_j$. This yields
the asserted result.
\end{proof}
\bigskip

We remark that the above argument also allows us to handle the
case when we replace the left side of (\ref{E:HKLRT}) by
$\prod_{i=1}^n y(z+c_i)$, which gives the same conclusion that
$\max\{p,\,q\}\le n$. This case was also considered in  \cite{AHH}
and \cite{HKLRT}.

We now consider the growth of meromorphic solutions to general
linear difference equation (\ref{E:higher-order-linear}).

\begin{theorem}\label{T:appl2}
Let $A_0(z),\, \cdots A_n(z)$ be entire functions such that there
exists an integer $\ell,\, 0\le \ell\le n$, such that
\begin{equation}
\label{E:appl2-1}
  \sigma (A_\ell)> \max_{\substack{0\le \ell\le n\\ j\not=\ell}}\{\sigma(A_j)\}.
\end{equation}
If $f(z)$ is a meromorphic solution to
\begin{equation}
\label{E:appl2-higher-order-linear}
    A_n(z)y(z+n)+\cdots+A_{1}(z)y(z+1) +A_0y(z)=0,
\end{equation}
then we have $\sigma(f)\ge \sigma(A_\ell)+1$.
\end{theorem}

\begin{proof} Let us choose $\sigma$ in relation to (\ref{E:appl2-1})
 so that
\begin{equation}
\label{E:proof-appl2-0}
    \max_{\substack{0\le \ell\le n\\
    \ell\not=j}}\{\sigma(A_i)\}<\sigma<\sigma (A_\ell)
\end{equation}
holds. Let us suppose that $f(z)$ is a finite order meromorphic
solution to (\ref{E:appl2-higher-order-linear}) such that
\begin{equation}
\label{E:proof-appl2-01}
    \sigma(f)<\sigma(A_\ell)+1.
\end{equation}
We divide through the equation (\ref{E:appl2-higher-order-linear})
by $f(z+\ell)$ to get
\begin{equation}
\label{E:proof-appl2-1}
    A_n(z)\frac{f(z+n)}{f(z+\ell)}+\cdots
    +A_\ell(z)+\cdots+A_0(z)\frac{f(z)}{f(z+\ell)}=0.
\end{equation}
Since (\ref{E:proof-appl2-0}) and (\ref{E:proof-appl2-01}) hold,
so we may choose $\varepsilon>0$ such that the inequalities
\begin{equation}
\label{E:proof-appl2-1.5}
   \sigma(f)+2\varepsilon<\sigma(A_\ell)+1\quad \text{and}\quad
    \sigma+2\varepsilon<\sigma (A_\ell),
\end{equation}
hold simultaneously. With the $\varepsilon>0$ as given in
(\ref{E:proof-appl2-1.5}), then (\ref{E:discrete-quotient-2})
gives, when $0\le j<\ell$ or $\ell< j\le n$,
\begin{equation} \label{E:proof-appl2-2}
    m\Big(r, \frac{f(z+j)}{f(z+\ell)}\Big)\le
    O\big(r^{\sigma(f)-1+\varepsilon}\big).
\end{equation}

Then we deduce from (\ref{E:discrete-proximity}),
(\ref{E:proof-appl2-2}) and (\ref{E:proof-appl2-1.5}) that
\begin{equation}
\label{E:proof-appl2-4}
\begin{split}
    m(r,\,A_\ell) &\le \sum_{\substack{0\le j\le n,\\ j\not= \ell}}
    m\Big(r, \frac{f(z+j)}{f(z+\ell)}\Big)+
    \sum_{j\not= \ell} m(r, A_j)\\
    &\le
    O\big(r^{\sigma(f)-1+\varepsilon}\big)+O\big(r^{\sigma+\varepsilon}\big)\\
    &\le o\big(r^{\sigma(A_\ell)-\varepsilon}\big).
\end{split}
\end{equation}
A contradiction.
\end{proof}

We next show how to use the Theorem \ref{T:appl2} to settle a
problem of Whittaker \cite{Whit} concerning linear difference
equations.

\begin{corollary}
\label{C:appl-1} Let $\sigma$ be a real number, and let $\Psi(z)$
be a given entire function with with order $\sigma(\Psi)=\sigma$.
Then the equation
\begin{equation}
\label{E:1st-order-2}
    F(z+\eta)=\Psi(z) F(z),
\end{equation}
admits a meromorphic solution of order $\sigma(F)=\sigma+1$.
\end{corollary}

\begin{proof} Whittaker \cite[\S 6]{Whit} constructed a meromorphic solution
$F(z)$ to the equation (\ref{E:1st-order-2}) and the solution has
order $\sigma(F)\le \sigma(\Psi)+1$. Since $\Psi$ is entire, and
it certainly satisfies the assumption (\ref{E:appl2-1}) and this
leads to the conclusion that $\sigma(F)\ge \sigma(\Psi)+1$. This
completes the proof.
\end{proof}

\begin{theorem}\label{T:appl3}
Let $P_0(z),\, \cdots P_n(z)$ be polynomials such that there
exists an integer $\ell,\, 0\le \ell\le n$ so that
\begin{equation}
\label{E:appl3-1}
  \deg (P_\ell)> \max_{\substack{0\le \ell\le n\\ j\not=\ell}}\{\deg(P_j)\}
\end{equation}
holds. Suppose $f(z)$ is a meromorphic solution to
\begin{equation}
\label{E:appl3-higher-order-linear}
    P_n(z)y(z+n)+\cdots+P_{1}(z)y(z+1) +P_0y(z)=0,
\end{equation}
then we have $\sigma(f)\ge 1$.
\end{theorem}
\medskip

\begin{proof} We assume that the equation
(\ref{E:appl3-higher-order-linear}) admits a meromorphic solution
$f(z)$ with $\sigma(f)<1$. We now divide through the difference
equation (\ref{E:appl3-higher-order-linear}) by $f(z+\ell)$ to
obtain
\begin{equation}
\label{E:proof-appl3-1}
    P_n(z)\frac{f(z+n)}{f(z+\ell)}+\cdots
    +P_\ell(z)+\cdots+P_0(z)\frac{f(z)}{f(z+\ell)}=0.
\end{equation}
We note that since $\sigma(f)<1$, so let us choose an
$\varepsilon>0$ so that $\varepsilon<1-\sigma(f)$, and Corollary
\ref{C:pointwise-1} implies that both when $0\le j< \ell$ or
$\ell< j\le n$ hold, then
\begin{equation}
\label{E:proof-appl3-2}
    \bigg|\frac{f(z+j)}{f(z+\ell)}\bigg|\leq
    \exp\big(r^{\sigma-1+\varepsilon}\big)=\exp\big({o(1)}\big)
\end{equation}
also holds outside a possible set of $r$ of finite logarithmic
measure. We deduce that (\ref{E:proof-appl3-2}) is bounded outside
a possible set of $r$ of finite logarithmic measure.

We now apply the (\ref{E:proof-appl3-2}) to equation
(\ref{E:proof-appl3-1}) and this gives
\begin{equation}
\label{E:proof-appl3-4}
    %\begin{split}
    |P_\ell(z)| \le \sum_{\substack{0\le j\le n,\\ j\not=\ell}}
    |P_j(z)|\bigg|\frac{f(z+j)}{f(z+\ell)}\bigg|
    %+\sum_{0\le  j<\ell}\|P_j(z)\|\bigg|\frac{f(z+j)}{f(z+\ell)}\bigg|\\
    \le O(1)\sum_{\substack{0\le j\le n,\\ j\not=\ell}}
    |P_j(z)|,
    %\end{split}
\end{equation}
as $|z|\to \infty$, outside a possible set of $r$ of finite
logarithmic measure. A contradiction to the assumption
(\ref{E:appl3-1}).
\end{proof}

We consider the following examples showing the sharpness of the
above theorems.

\begin{example} Ruijsenaars \cite{Ruij} considers the equation
\begin{equation}
\label{E:Ruij-1}
     F(z+ia/2)=\Phi(z)F(z-ia/2),
\end{equation}
where  $a>0$, $\Phi(z)=2\cosh {\pi z/b}$ and $b>0$. The solution
\begin{equation}
\label{E:Ruij-2}
    G_{\text{hyp}}(a,\,b;
    z)=\exp\left(i\Big(\int_0^\infty
    \frac{\sin(2yz)}{2\sinh(ay)\sinh(by)}-\frac{a}{aby}\Big)\,\frac{dy}{y}\right),\quad
    |\Im z|<(a+b)/2,
\end{equation}
which has no zeros and poles in $ |\Im z|<(a+b)/2$, can be
continued meromorphically to the whole complex plane via the
equation (\ref{E:Ruij-1}). The poles and zeros of (\ref{E:Ruij-2})
are given, respectively, by
\begin{equation}
\label{E:Ruij-3}
    z=-i(k+1/2)a-i(\ell+1/2)b,\quad
    z=-i(k+1/2)a+i(\ell+1/2)b,\quad k,\, \ell\in \mathbf{N}.
\end{equation}
The function (\ref{E:Ruij-2}) is called the \textit{hyperbolic
gamma function}. It follows from (\ref{E:Ruij-3}) that the order
of $G_{\text{hyp}}(r,\,a;\,z)$ is $2$. Thus we have
$\sigma(G_{\text{hyp}})=\sigma(\Phi)+1$. We would like to mention
that Ruijsenaars \cite{Ruij} also considers the equation
(\ref{E:Ruij-1}) where $\Phi$ is the \textit{trigonometric gamma
function} and the \textit{elliptic gamma function} respectively.
We again have $\sigma(G_{\text{ell}})=\sigma(\Phi)+1$, and
$\sigma(G_{\text{trig}})=\sigma(\Phi)+1$ to hold. We also remark
that all the three types of generalized gamma functions mentioned
above converge to the Euler Gamma function while taking suitable
limits of the parameters.
\end{example}

The following equation was considered in Hayman and Thatcher
\cite{Hay90}, which has a different form from (\ref{E:Ruij-1}).
Let $H>0$, then the equation
\begin{example}
\begin{equation}
\label{E:Hayman-1}
    F(z)=(1+H^z)F(z+1),
\end{equation}
admits a meromorphic solution of the form
\begin{equation}
\label{E:Hayman-2}
    F_1(z)=\prod_{n=1}^\infty \big(1+H^{z-n}\big)^{-1}
\end{equation}
with simple poles at
\begin{equation}
\label{E:Hayman-3}
    z_{k,\,n}=n+\frac{(2k+1)\pi i}{\log H},
\end{equation}
where $n=1,\, 2,\, 3,\,\cdots$ and $k$ is an integer \cite[Thm.
1]{Hay90}. It follows from (\ref{E:Hayman-3}) that $F_1$ has order
2, giving $\sigma(F_1)=\sigma (1+H^z)+1$ so that the ``equality"
holds in Theorem \ref{T:appl2} again.
\end{example}

The next example shows that the assumption (\ref{E:appl3-1}) where
only one coefficient is allowed to have the highest degree is the
best possible.
\begin{example}[\cite{Is-Ya}] Let
\begin{equation}
\label{E:eg3-1}
    \Delta^n f(z)=\sum_{j=0}^n\binom{n}{j} (-1)^{n-j}f(z+j)
\end{equation}
and hence
\begin{equation}
\label{E:eg3-2}
    f(z+n)=\sum_{j=0}^n\binom{n}{j} \Delta^j f(z).
\end{equation}
Then the equation
\begin{equation}
    z(z-1)(z-2)\Delta^3f(z-3)+z(z-1)\Delta^2 f(z-2)+z\Delta
    f(z-1)+(z+1)f(z)=0
\end{equation}
admits an entire solution of order 1/3. In fact, it is shown in
\cite{Is-Ya} that
\[
    \log M(r,\, f)=L r^{1/3}\big((1+o(1)\big).
\]
By making use the relation (\ref{E:eg3-2}), we can rewrite the
equation (\ref{E:eg3-1}) to an equation of the form
(\ref{E:appl3-higher-order-linear}) with
\[
\deg P_3=\deg P_2=\deg P_1=\deg P_0=3.
\]
Thus, there are more than one polynomial coefficients having the
same degree ($>0$) and the equation admits an entire solution of
order $<1$. The above example shows that we cannot drop the
assumption (\ref{E:appl3-1}) in Theorem \ref{T:appl3}.
\end{example}

We finally remark that the equation (\ref{E:1st-order}) and its
solution (\ref{E:gamma-generalized}) show that the lower order one
estimate in the Theorem \ref{T:appl3} is again the best possible.

\section{Discussion}

In this paper, we have discussed in detail some basic properties
of $T\big(r,\,f(z+\eta)\big)$, for a fixed $\eta$. In particular,
we have shown in Theorem \ref{T:char} that the relation
(\ref{E:char-asy}) holds for finite order meromorphic functions
and that the Theorem \ref{T:counting-example} shows that no such
relation (\ref{E:char-asy}) can hold for infinite order
meromorphic functions. The proof of (\ref{E:char-asy}) depends on
Theorem \ref{T:discrete-quotient} which can be viewed as a
discrete analogue of the classical logarithmic derivative estimate
given by Nevanlinna (\cite{Nev}, \cite{Hay64}) and the relation
(\ref{E:counting-est}) in Theorem \ref{T:counting} on the counting
function. These special properties of finite order meromorphic
functions distinguish themselves from general meromorphic
functions. They are in strong agreement with the integrability
detector of difference equations proposed in \cite{AHH}.

It is worthwhile to note that the integrability test by Nevanlinna
theory proposed in \cite{AHH} is being complex analytic in nature,
which is in stark contrast when compared to several other major
integrability tests for difference equations (\cite{CM},
\cite{GRP}, \cite{RGH}, \cite{RGTT}) proposed in the last decade.
In fact, the \textit{Nevanlinna test} seems more natural when
compared to the well-known complex analytic \textit{Painlev\'e
test} as an integrability test for second order ordinary
differential equations; see \cite[p. 362]{AC}. We mention that the
prime integrable difference equations are the discrete Painlev\'e
equations which can be obtained from the classical Painlev\'e
differential equations \cite{GLS} via suitable discretizations.

Although the investigation in \cite{AHH} is for non-linear second
order difference equations, we have found that it natural to
consider \textit{linear difference equations}. This is based on
the following facts. First our investigation leads us to give an
answer of a Whittaker's problem (Corollary \ref{C:appl-1}), which
is amongst the most basic results of first order difference
equations from the viewpoint of Nevanlinna theory. Second, we use
our main result (Theorem \ref{T:char}) to give a simple proof of
the main result (Thm. \ref{T:HKLRT}) in \cite{AHH}.

Linear difference equations are generally accepted as integrable.
From the viewpoint in \cite{AHH}, it is therefore natural to
demand that meromorphic solutions to linear equations should also
be of finite order of growth. However, the discussion in \S1 and
Theorem \ref{T:appl2} indicate that meromorphic solutions to
(\ref{E:1st-order}) could have an arbitrarily fast growth. We give
a lower bound order estimate of a finite order meromorphic
solution, if any, of a linear equation (Theorem \ref{T:appl2}; see
also Theorem \ref{T:appl3}). Thus one must impose certain
\textit{minimal growth} condition to single out the minimal
solution (and finite order, if any). The question here is that
what determines a \textit{minimal solution}. The Whittaker theorem
(Corollary \ref{C:appl-1}) shows that minimal solution always
exist for first order equation with an arbitrary entire
coefficient in terms of order of growth. The problem of
\textit{minimal solution} is investigated in \cite{CR} for certain
first order difference equations where the meromorphic solutions
has prior growth restriction in an infinite strip. The distinction
of different \textit{minimal solutions} is also discussed in
\cite{Hay90}.

\vskip1cm \noindent\textbf{Acknowledgements} The authors would
like to thank Dr. Patrick Ng of the University of Hong Kong who
brought to the attention of the authors of the preprint by R. G.
Halburd and R. J. Korhonen \cite{HK-1}. The authors would also
like to thank Dr. Mourad Ismail and the referee for useful
comments to our paper.

\bibliographystyle{amsplain}

%\vspace{-1.3cm}
\end{document}